\documentclass{article}
\usepackage{amsmath,amsthm}
\usepackage{graphicx,epstopdf,epsfig,multirow,epic,bm}
\usepackage{subfigure}

\normalsize \rm
\begin{document}

\begin{center}
{\Large  \textbf {Exact mean first-passage time on generalized Vicsek fractal}}\\[12pt]
{\large \quad Fei Ma$^{a,}$\footnote{~The author's E-mail: mafei123987@163.com. } \;  Xiaomin Wang$^{a}$ \; Ping Wang$^{b,c,d}$ \;and \; Xudong Luo$^{e}$}\\[6pt]
{\footnotesize $^{a}$ School of Electronics Engineering and Computer Science, Peking University, Beijing 100871, China\\
$^{b}$ National Engineering Research Center for Software Engineering, Peking University, Beijing, China\\
$^{c}$ School of Software and Microelectronics, Peking University, Beijing  102600, China\\
$^{d}$ Key Laboratory of High Confidence Software Technologies (PKU), Ministry of Education, Beijing, China\\
$^{e}$ College of Mathematics and Statistics, Northwest Normal University, Lanzhou  730070, China}\\[12pt]
\end{center}

\begin{quote}
\textbf{Abstract:}  Fractal phenomena may be widely observed in a great number of complex systems. In this paper, we revisit the well-known Vicsek fractal, and study some of its structural properties for purpose of understanding how the underlying topology influences its dynamic behaviors. For instance, we analytically determine the exact solution to mean first-passage time for random walks on Vicsek fractal in a more light mapping-based manner than previous other methods, including typical spectral technique. More importantly, our method can be quite efficient to precisely calculate the solutions to mean first-passage time on all generalized versions of Vicsek fractal generated based on an arbitrary allowed seed, while other previous methods suitable for typical Vicsek fractal will become prohibitively complicated and even fail. Lastly, this analytic results suggest that the scaling relation between mean first-passage time and vertex number in generalized versions of Vicsek fractal keeps unchanged in the large graph size limit no matter what seed is selected. \\

\textbf{Keywords:}  Vicsek fractal, Random walks, Mean first-passage time. \\

\end{quote}

\vskip 1cm

\section{Introduction}

Fractals are intriguing and may be popularly seen in a great deal of complex systems \cite{Mandlebrot-1982}. There is a long history of investigation of various kinds of fractal models. The best known are Koch Curve \cite{Paramanathan-2010}, Appollonian networked model \cite{Andrade-2005}, Sierpinski carpet \cite{Barlow-1999}, T-graph \cite{Agliari-2008} as well as Vicsek fractal \cite{Vicsek-1983}-\cite{Jayanthi-1992}, and so forth. The early related works have uncovered some structural properties, for instance, fractal dimension and spectral dimension, on these fractal models mentioned above. Nonetheless, in this paper, we intend to consider the famous Vicsek fractal, a classic family of regular hyperbranched polymers that have a great number of diverse applications in this field of polymer physics \cite{Gurtovenko-2005}, and study its some topological parameters, including computation of mean first-passage time for random walks on Vicsek fractal, in order to understand how the underlying topologies of models of this kind influence their dynamic behaviors.

Random walks, as a classic representative featuring diffusion phenomena on complex networks, have  attracted an increasing attention from a wide range of scientific fields \cite{Noh-2004}-\cite{Masuda-2017}, including statistic physics, applied mathematics, theoretical chemistry and theoretical computer science, and so on. The most significant reason for this is that there are a great number of practical applications closely correlated to random walks, such as, protein folding in cell, the motions of molecules in rarified gases, information flow in social networks, the behaviour of stochastic search algorithms in the World Wide Web (WWW), to name just a few \cite{Perkins-2014}. In the language of mathematics, a walker doing random walks on a given network will hop on each of its neighbor vertex set from the current location $u$ with the probability $1/k_{u}$ where $k_{u}$ is the degree of vertex $u$. In general, such a dynamic process can be described by the transition matrix $\mathbf{P}=\mathbf{D}^{-1}\mathbf{A}$ of network under consideration where $\mathbf{A}$ is the adjacency matrix and $\mathbf{D}$ represents the corresponding diagonal matrix (defined in detail later). The key issue in studying random walks on networks is to estimate a topological parameter, called mean first-passage time. It is well known that such a parameter has given rise to a growing number of theoretical investigations over the past several decades \cite{Hwang-2012}-\cite{Condamin-2007}. Noticeably, throughout this paper, all graphs (models) considered are simple and connected, and the terms graph and network are used indistinctly. Symbol $[a,b]$ is a set consisting of integers $\{a,a+1,\cdots,b-1,b\}$.

\begin{figure*}
\centering
\subfigure[]{
\begin{minipage}[t]{0.5\linewidth}
\centering
\includegraphics[width=7cm]{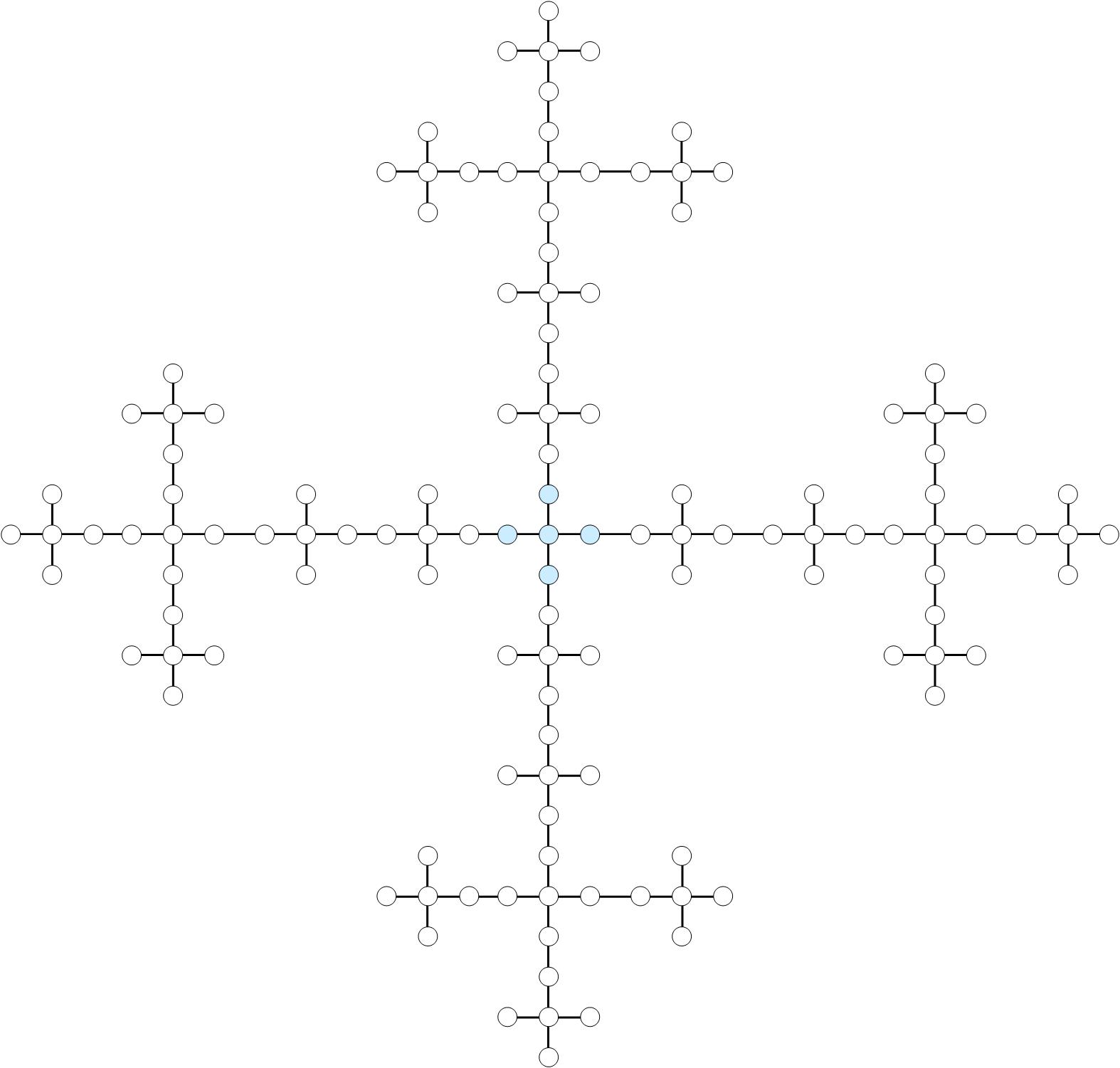}
%\caption{fig1}
\end{minipage}%
}%
\subfigure[]{
\begin{minipage}[t]{0.5\linewidth}
\centering
\includegraphics[width=7cm]{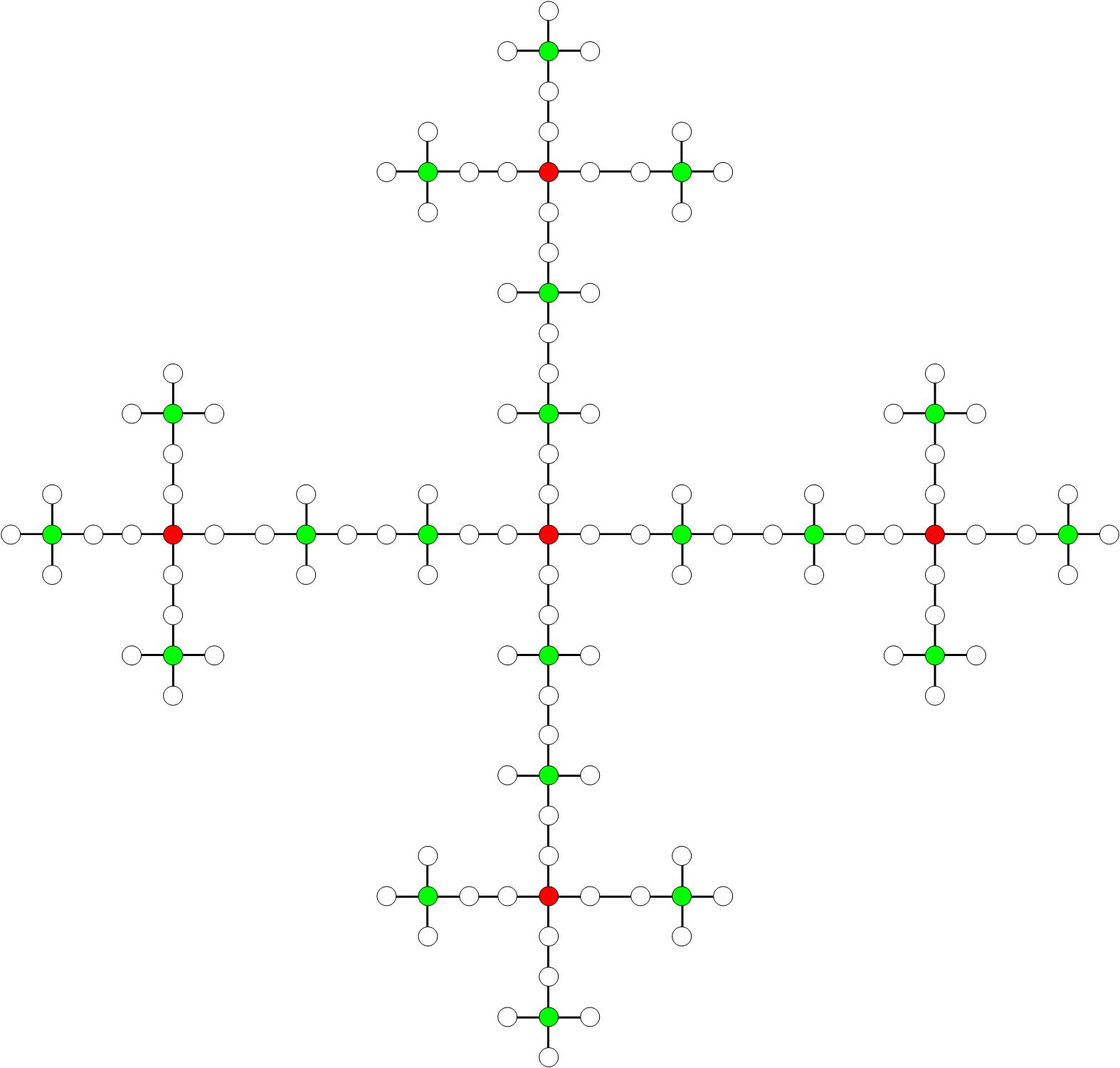}
%\caption{fig2}
\end{minipage}
}%
\caption{\label{fig:wide} (Color online) The diagram of Vicsek fractal $\mathrm{V}_{2,4}$. The filled vertices in panel (a) denote the initial structure of Vicsek fractal $\mathrm{V}_{0,4}$. According to Vicsek fractal operation $V_{4}$ illustrated in Fig,2, the panel (b) shows a detailed construction of Vicsek fractal $\mathrm{V}_{2,4}$ where the red vertices represent the initial vertices at time step $t=0$, followed by green vertices indicating ones created at time step $t=1$, and the left ones are introduced at time step $t=2$.}
\end{figure*}

In fact, concerted efforts have been devoted to discussing mean first-passage time on networks including fractal models in order to unveil the scaling law of quantity with the vertex number of networks. The results in recent paper \cite{Perkins-2014} confirm that there are precisely three mutually exclusive possibilities for the form of the path distribution: finite, stretched exponential and power law, found in any random walk on a finite network. In particular, some deterministic networked models are also considered. Among which, the T-graph is in depth discussed and thus its corresponding solution to mean first-passage time has been reported in \cite{Agliari-2008}. The most close to our interest is work done by Wu \emph{et al} in \cite{Wu-2012} where the typical Vicsek fractal is chose as an object. As a result, the corresponding solution to mean first-passage time has been obtained based on its own self-similarity. As will be demonstrated shortly, we consider a more general models covering typical Vicsek fractal and then determine the formula for mean first-passage time on such generalized versions built based on an arbitrary allowed seed in a more convenient manner developed later. However, many other previous methods suitable for typical Vicsek fractal will become prohibitively complicated and even fail to address this issue. The main purpose of this paper is to fill this gap. Finally, the results firmly suggest that the scaling relation between mean first-passage time and vertex number in generalized versions of such type remains unchanged in the large graph size limit no matter what that seed is selected.

The rest of this paper can be organized as follows. In Section 2, we discuss random walks on Vicsek fractal. First of all, some fundamental notations are introduced. And then, we propose three algorithms for constructing Vicsek fractal. The first two of them are designed for producing typical Vicsek fractal. The left one is for establishing the generalized version of typical Vicsek fractal. Next, we analytically calculate the closed-form solution to mean first-passage time for random walks on generalized version of typical Vicsek fractal in a mapping-based manner. Particularly, the corresponding solution on typical Vicsek fractal is easily covered by our results. In addition, we analyze the scaling relations on some quantities. Finally, we draw the conclusion in Section 3.

\section{Random walks on Vicsek fractal }

In this section, we are going to study random walks on Vicsek fractal. For our purpose, the first subsection aims at introducing some basic notations used later, including the matrix representation of graph. The following subsection will propose three algorithms for establishing Vicsek fractal. Lastly, we determine exact solutions to some quantities that we are interesting in.

\subsection{Basic notations}

\subsubsection{Matrix representation of graph}
Consider a graph $G(V,E)$, it is often described by the adjacency matrix $\mathbf{A}=(a_{uv})$ in the next form

$$a_{uv}=\left\{\begin{aligned}&1, \quad\text{vertex $u$ is adjacent to $v$;}\\
&0,\quad\text{otherwise}.
\end{aligned}\right.
$$
Based on such a description, some basic information related to a graph itself can be handily obtained. For instance, the degree $k_{u}$ of vertex $u$ is equal to $k_{u}=\sum_{v=1}^{|V|}a_{uv}$. And, the diagonal matrix, denoted by $\mathbf{D}$, may be defined as follows: the $i$th diagonal entry is $k_{i}$, while all non-diagonal elements are zero, i.e., $\mathbf{D}=\text{diag}[k_{1},k_{2},\dots,k_{|V|}]$.

\subsubsection{Wiener index}
For a given graph $G(V,E)$, the distance $d_{uv}$ between a pair of vertices $u$ and $v$, also defined as shortest path length, is the edge number of any shortest path joining vertices $u$ and $v$. The sum over distances of all possible vertex pairs, commonly called \emph{Wiener index} $\mathcal{W}$ \cite{Georgakopoulos-2017}, is by definition given by

\begin{equation}\label{eqa:MF-2-A-1}
\mathcal{W}=\sum_{u,v;u\neq v}d_{uv}=\frac{1}{2}\sum_{u\in V}\sum_{v\in V}d_{uv}.
\end{equation}
Thus, the average path length $\langle\mathcal{W}\rangle$ is viewed as
\begin{equation}\label{eqa:MF-2-A-2}
\langle\mathcal{W}\rangle=\frac{\mathcal{W}}{|V|(|V|-1)/2}
\end{equation}
in which we denote by $|V|$ vertex number of graph $G(V,E)$.

\subsubsection{Mean first-passage time}

Similarly, the crucial issue in studying random walks on graph $G(V,E)$ is to measure a quantity, called \emph{mean first-passage time} $\mathcal{A}$. The definition of $\mathcal{A}$ is

\begin{equation}\label{eqa:MF-2-A-3}
\mathcal{A}=\frac{\sum_{u,v;u\neq v}F_{u\rightarrow v}}{|V|(|V|-1)}
\end{equation}
where $F_{u\rightarrow v}$ indicates \emph{first-passage time} took by a walker starting out from source vertex $u$ to first reach destination vertex $v$.

\subsection{Construction of Vicsek fractal}

Here, we first introduce the detailed construction of typical Vicsek fractal. Accordingly, the generalized versions will be built based on an arbitrary seed.

The typical Vicsek fractal, referred to as $\mathrm{V}_{t,s}$ ($s\geq2$), is constructed in an iterative manner as follows.

\textbf{\emph{Algorithm 1}}

At $t=0$, the seed,  denoted by $\mathrm{V}_{0,s}$, is a star having $s$ leaf vertices arranged in a cross-wise pattern.

At $t\geq1$, the next generation $\mathrm{V}_{t,s}$ is obtained from preceding model $\mathrm{V}_{t-1,s}$ in the following way. First of all, in order to obtain $\mathrm{V}_{1,s}$, one should creates $s$ duplications $\mathrm{V}^{i}_{0,s}$ ($i\in[1,s]$) of $\mathrm{V}_{0,s}$ that are successively placed around the periphery of the original $\mathrm{V}_{0,s}$. And then, one takes $s$ new edges to connect $\mathrm{V}_{0,s}$ to $\mathrm{V}^{i}_{0,s}$ in order. Such a procedure can be executed $t$ time steps until the desirable Vicsek fractal $\mathrm{V}_{t,s}$ is generated. As an illustrative example of Vicsek fractal $\mathrm{V}_{t,s}$, panel (a) in Fig.1 shows $\mathrm{V}_{2,4}$. $\hfill\qedsymbol$

In view of the representation in Algorithm 1, after $t$ time steps, the total number $|\mathrm{V}_{t,s}|$ of vertices of Vicsek fractal $\mathrm{V}_{t,s}$ is calculated in a recursive manner and follows

\begin{equation}\label{eqa:MF-2-B-1}
|\mathrm{V}_{t,s}|=(s+1)^{t+1}.
\end{equation}

In essence, the above-proposed Vicsek fractal $\mathrm{V}_{t,s}$ is also generated in another iterative fashion based on \emph{Vicsek fractal operation }$V_{s}$ defined below.

As illustrated in Fig.2, applying Vicsek fractal operation $V_{s}$ to a given edge $uv$ is to first insert two new vertices into edge $uv$, which makes edge $uv$ saturated, and to then connect $s-1$ new vertices to each end-vertex in order to make vertex degree equal to $s$ since each end-vertex has degree $1$ previously. In this situation, vertices $u$ and $v$ are considered $s$-saturated. With the help of Vicsek fractal operation $V_{s}$, we propose the next algorithm.

\begin{figure}
\centering
\includegraphics[height=6cm]{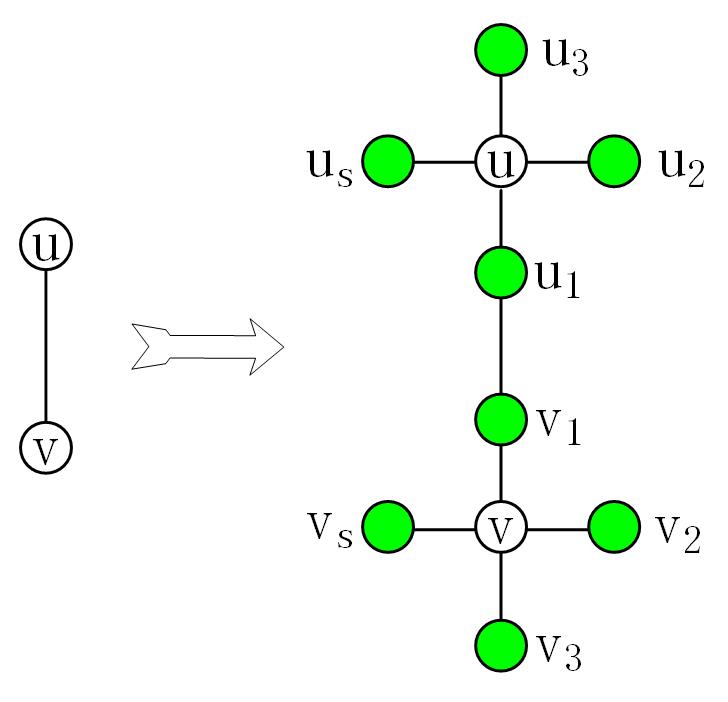}% Here is how to import EPS art
\caption{\label{fig:wide} (Color online)  The diagram of Vicsek fractal operation $V_{s}$ on an edge $uv$. The edge $uv$ in the left-hand side of arrow will result in a graph in the right-hand side of arrow using Vicsek fractal operation $V_{s}$.  }
\end{figure}

\textbf{\emph{Algorithm 2}}

As previously, the seed $\mathrm{V}_{0,s}$ is still a star having $s$ leaf vertices. At $t\geq1$, the next generation $\mathrm{V}_{t,s}$ can be obtained from previous model $\mathrm{V}_{t-1,s}$ by employing Vicsek fractal operation $V_{s}$ to each edge. Specifically, we first insert two new vertices into each edge in order to make all edges saturated. And then, we should connect $s-k_{v}$ new vertices to each vertex $v$ with degree $k_{v}$ such that vertex $v$ becomes $s$-saturated. Such a manipulation may be iteratively conducted $t$ time steps until the final Vicsek fractal $\mathrm{V}_{t,s}$ is produced. The panel (b) in Fig.1 plots an illustrative example $\mathrm{V}_{2,4}$. $\hfill\qedsymbol$

As we will show in the following, it is the construction described in Algorithm 2 that enables us to easily accomplish what we want.

Although the above both algorithms have successfully produced the typical Vicsek fractal $\mathrm{V}_{t,s}$, our goal is to study the generalized version $\mathcal{V}_{t,s}$ of Vicsek fractal $\mathrm{V}_{t,s}$. Accordingly, the previous results are completely covered by our ones. To this end, we must extend the typical Vicsek fractal $\mathrm{V}_{t,s}$. First, let us recall the concrete development of Vicsek fractal $\mathrm{V}_{t,s}$ for purpose of finding out some useful structural properties. The most striking is that there are only three types of vertices in Vicsek fractal $\mathrm{V}_{t,s}$ with respect to vertex degree. That is to say, degree value equal to $1,2$, and $s$ can be observed in Vicsek fractal $\mathrm{V}_{t,s}$ and others are absent.  Inspired by this, a natural generalization is obtained by Algorithm 3 as below.

\textbf{\emph{Algorithm 3}}

Let the seed be a tree whose vertices have degree values only equal to either $1, 2$, or $s$. The remainder is to manipulate  Vicsek fractal operation $V_{s}$ in a similar manner to that mentioned Algorithm 2 until the generalized version $\mathcal{V}_{t,s}$ is yielded. $\hfill\qedsymbol$

Note that the requirement about seed stated in Algorithm 3 can be able to further loose such that an arbitrary tree is allowed to choose as a seed. That is to say, for a given tree whose vertices have maximum degree $s$, one has the ability to produce generalized version $\mathcal{V}_{t,l}$ ($l\geq s$) utilizing Vicsek fractal operation $V_{l}$ $t$ time steps. For illustration purpose, an example $\mathcal{V}_{1,4}$ is plotted in Fig.3. It is straightforward to see that after employing Vicsek fractal operation $V_{l}$ only once, all vertices in the seed have degree $l$ now. As a result, we focus particularly on generalized version $\mathcal{V}_{t,s}$, and derive some solutions to structural parameters including mean first-passage time.

\begin{figure}
\centering
\includegraphics[height=5cm]{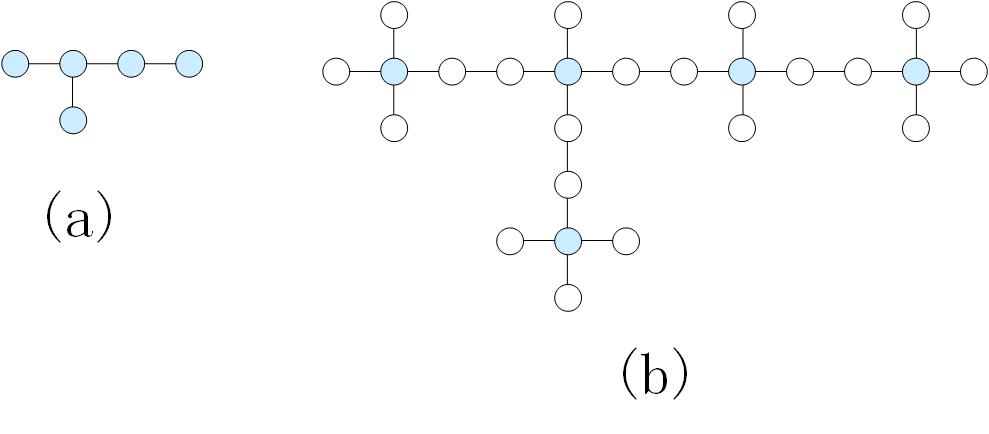}% Here is how to import EPS art
\caption{\label{fig:wide} (Color online)  The diagram of generalized version $\mathcal{V}_{t,s}$. The panel (a) shows a tree whose maximum vertex degree is equal to $3$ that is selected as a seed. With Vicsek fractal operation $V_{4}$, the resulting $\mathcal{V}_{1,4}$ is a model as shown in panel (b).    }
\end{figure}

\subsection{Exact solution to mean first-passage time $\mathcal{A}_{\mathcal{V}_{t,s}}$ for random walks on Vicsek fractal $\mathcal{V}_{t,s}$}

Clearly, all Vicsek fractals proposed including typical factal $\mathrm{V}_{t,s}$ share tree structure. In what follows, we introduce a helpful lemma on tree.

\textbf{Lemma 1} For a tree $T$, there is in fact a connection between Wiener index $\mathcal{W}$ and mean first-passage time $\mathcal{A}$, namely,

\begin{equation}\label{eqa:MF-2-C-1}
\mathcal{A}=\frac{2\mathcal{W}}{|T|}
\end{equation}
where $|T|$ is vertex number of tree itself \cite{Ma-2020}.

This leads to a fact that determining exact solution to mean first-passage time for random walks on Vicsek fractal $\mathcal{V}_{t,s}$ is to equivalently measure the corresponding Wiener index. As known, such a helpful demonstration has been widely used to obtain some closed-form solutions to mean first-passage time on many intriguing trees including T-graph. As shown in the following, we also attempt to follow such a thought to precisely derive solution to mean first-passage time $\mathcal{A}_{\mathcal{V}_{t,s}}$ of Vicsek fractal $\mathcal{V}_{t,s}$.

To answer this issue mentioned above, let us first consider a simple case, i.e., looking for a closed-form expression of Wiener index $\mathcal{W}_{\mathcal{V}_{1,s}}$ on Vicsek fractal $\mathcal{V}_{1,s}$ in terms of the known Wiener index $\mathcal{W}$ of its seed (i.e., a tree $T$) with $n$ vertices. Beginning with proceeding this issue, we have to make use of a trivial result that there is a unique path between arbitrary pair of vertices in tree. Indeed, this characteristic allows us to easily calculate what we want.

\textbf{Theorem 2} The Wiener index of generalized version $\mathcal{V}_{1,s}$ is

\begin{equation}\label{eqa:MF-2-C-7}
\mathcal{W}_{\mathcal{V}_{1,s}}=3(s+1)^{2}\mathcal{W}+(s^{2}-s-2)n^{2}+(s+2)n.
\end{equation}

\emph{\textbf{Proof}} Without loss of generality, a distinct label $v_{i}$ ($i\in[1,n]$) is assigned to each vertex in the seed $T$ (i.e., $\mathcal{V}_{0,s}$) such that we can take advantage of symbol $P_{v_{i}v_{j}}$ to represent that unique path joining vertex $v_{i}$ to vertex $v_{j}$. Applying Vicsek fractal operation $V_{s}$ to model $\mathcal{V}_{0,s}$ leads to a family of new vertices. Among which, as shown in Fig.2, each vertex $v_{i}$ will be connected to $s$ new vertices that are in turn labeled by using $v_{i}^{j}$ ($j\in[1,s]$). By far, we are ready to deal with computation of quantity $\mathcal{W}_{\mathcal{V}_{1,s}}$. The complete computation consists of five stages as shown in the following.

\emph{Case 1} For a given couple of vertices $v_{i}$ and $v_{j}$, it is not hard to see that its distance $d'_{v_{i}v_{j}}$ in generalized version $\mathcal{V}_{1,s}$ is three times larger than distance $d_{v_{i}v_{j}}$ in seed $\mathcal{V}_{0,s}$, suggesting

\begin{equation}\label{eqa:MF-2-C-2}
\mathcal{W}_{1}=\sum_{v_{i},v_{j}}d'_{v_{i}v_{j}}=3\mathcal{W}.
\end{equation}

\emph{Case 2} From illustration shown in Fig.2, there are $s$ new vertices $v_{i}^{j}$ connected to each vertex $v_{i}$. Intuitively, the underlying structure induced by such $s+1$ vertices is star-like. So, the sum $\mathcal{W}_{2i}$ of distances between central vertex $v_{i}$ and each leaf vertex $v_{i}^{j}$ is equal to $s$. There will be $n$ star-like underlying subgraphs of such kind in generalized version $\mathcal{V}_{1,s}$, implying

\begin{equation}\label{eqa:MF-2-C-3}
\mathcal{W}_{2}=\sum_{i=1}^{n}\mathcal{W}_{2i}=sn.
\end{equation}

\emph{Case 3} Following the above discussion, we now consider a case of determining exact solution to sum $\mathcal{W}_{3}$ of distances $d'_{v_{i}v_{j}^{r}}$ over all possible vertex pairs $v_{i}$ and $v_{j}^{r}$. It is worth noting that subscript $i$ is distinct with subscript $j$. To put this another way, vertex $v_{i}$ is different from central vertex $v_{j}$ connected to vertex $v_{j}^{r}$. To address this issue, we need to take a mapping $f_{1}$ between paths $P_{v_{i}v_{j}^{r}}$ and $P_{v_{i}v_{j}}$. Obviously, there exist $s$ distinct paths $P_{v_{i}v_{j}^{r}}$ of this kind. Taking into account the underlying structure of generalized version $\mathcal{V}_{1,s}$, we construct mapping $f_{1}$ as follows

$$f_{1}:\; P_{v_{i}v_{j}^{r}}\mapsto \; P_{v_{i}v_{j}}.$$
The performances of mapping $f_{1}$ encompass that (1) There must be $(s-1)$ paths $P_{v_{i}v_{j}^{r}}$ degraded into path $P_{v_{i}v_{j}}$ by removing an additional edge $v_{j}v_{j}^{r}$ under mapping $f_{1}$, (2) The remaining that path $P_{v_{i}v_{j}^{r'}}$ can be extended into path $P_{v_{i}v_{j}}$ by adding an additional edge $v_{j}^{r'}v_{j}$ under mapping $f_{1}$. In fact, there is also a completely the same mapping $f_{1}$ between paths $P_{v_{i}^{r}v_{j}}$ and $P_{v_{i}v_{j}}$. Due to similar analysis, we omit detailed description. Taken together, we have

\begin{equation}\label{eqa:MF-2-C-4}
\begin{aligned}\mathcal{W}_{3}&=2(s-1)\left[\mathcal{W}_{1}+\frac{n(n-1)}{2}\right]+2\left[\mathcal{W}_{1}-\frac{n(n-1)}{2}\right]\\
&=2s\mathcal{W}_{1}+(s-2)n(n-1)
\end{aligned}.
\end{equation}

\emph{Case 4} As stated in case 2, there are $n$ star-like underlying subgraphs in generalized version $\mathcal{V}_{1,s}$. The next task is to measure sum $\mathcal{W}_{4i}$ of distances between two arbitrary leaf vertices in an identical star-like underlying subgraph. Armed with these discussions, we can write

\begin{equation}\label{eqa:MF-2-C-5}
\mathcal{W}_{4}=\sum_{i=1}^{n}\mathcal{W}_{4i}=s(s-1)n.
\end{equation}

\emph{Case 5} The left task is to measure distance $d'_{v_{i}^{l}v_{j}^{r}}$ of both vertices $v_{i}^{l}$ and $v_{j}^{r}$ where vertices $v_{i}^{l}$ and $v_{j}^{r}$ are from two distinct star-like underlying subgraphs. In other words, subscript $i$ is different from subscript $j$. Generalized version $\mathcal{V}_{1,s}$ has $s^{2}$ vertex pairs of this type for vertices $v_{i}$ and $v_{j}$ in total. Motivated by developing proof of quantity $\mathcal{W}_{3}$, we also take a mapping $f_{2}$ between paths $P_{v_{i}^{l}v_{j}^{r}}$ and $P_{v_{i}v_{j}}$, as follows

$$f_{2}:\; P_{v_{i}^{l}v_{j}^{r}}\mapsto \; P_{v_{i}v_{j}}.$$

Here, as above, the performances of mapping $f_{2}$ are in fact shown in four aspects.

\begin{itemize}

\item There must be $(s-1)^{2}$ paths $P_{v_{i}^{l}v_{j}^{r}}$ reduced into path $P_{v_{i}v_{j}}$ by removing two additional edges $v_{i}v_{i}^{l}$ and $v_{j}v_{j}^{r}$ under mapping $f_{2}$.

\item There must exist $(s-1)$ paths $P_{v_{i}^{l}v_{j}^{r'}}$ degraded to path $P_{v_{i}v_{j}}$ by both removing an additional edges $v_{i}v_{i}^{l}$ and adding an edge $v_{j}^{r'}v_{j}$ under mapping $f_{2}$.

\item There must be $(s-1)$ paths $P_{v_{i}^{l'}v_{j}^{r}}$ transformed into path $P_{v_{i}v_{j}}$ by both removing an additional edges $v_{j}v_{j}^{r}$ and adding an edge $v_{i}^{l'}v_{i}$ under mapping $f_{2}$.

\item The left unique path $P_{v_{i}^{l'}v_{j}^{r'}}$ may be naturally extended into path $P_{v_{i}v_{j}}$ through adding two additional edges $v_{i}^{l'}v_{i}$ and $v_{j}^{r'}v_{j}$ under mapping $f_{2}$.

\end{itemize}

All discussions upon four aspects above together yields

\begin{equation}\label{eqa:MF-2-C-6}
\begin{aligned}\mathcal{W}_{5}&=(s-1)^{2}\left[\mathcal{W}_{1}+2\times\frac{n(n-1)}{2}\right]+2(s-1)\mathcal{W}_{1}+\left[\mathcal{W}_{1}-2\times\frac{n(n-1)}{2}\right]\\
&=(s-1)^{2}[\mathcal{W}_{1}+n(n-1)]+2(n-1)\mathcal{W}_{1}+\mathcal{W}_{1}-n(n-1)
\end{aligned}.
\end{equation}

So far, we have exhaustively enumerated distances over all possible vertex pairs from generalized version $\mathcal{V}_{1,s}$. Therefore, the Wiener index $\mathcal{W}_{\mathcal{V}_{1,s}}$ is given by in the next form,

$$\mathcal{W}_{\mathcal{V}_{1,s}}=\sum_{i=1}^{5}\mathcal{W}_{i}.$$$\hfill\qedsymbol$

Accordingly, the iterative construction of generalized version $\mathcal{V}_{t,s}$ tells us the following equation

\begin{equation}\label{eqa:MF-2-C-8}
\mathcal{W}_{\mathcal{V}_{t,s}}=3(s+1)^{2}\mathcal{W}_{\mathcal{V}_{t-1,s}}+(s^{2}-s-2)|\mathcal{V}_{t-1,s}|^{2}+(s+2)|\mathcal{V}_{t-1,s}|
\end{equation}
where $|\mathcal{V}_{t,s}|$ is the total number of vertices in generalized version $\mathcal{V}_{1,s}$, and is calculated equal to

\begin{equation}\label{eqa:MF-2-C-9}
|\mathcal{V}_{t,s}|=n(s+1)^{t}.
\end{equation}

Combining Eqs.(\ref{eqa:MF-2-C-8}) and (\ref{eqa:MF-2-C-9}) yields the explicit solution to Wiener index of generalized version $\mathcal{V}_{t,s}$, as follows.

\textbf{Theorem 3} The Wiener index of generalized version $\mathcal{V}_{t,s}$ is

\begin{equation}\label{eqa:MF-2-C-10}
\mathcal{W}_{\mathcal{V}_{t,s}}=3^{t}(s+1)^{2t}\mathcal{W}+\frac{n^{2}(3^{t}-1)(s+1)^{2(t-1)}}{2}+\frac{n(s+1)^{t-1}[3^{t}(s+1)^{t}-1]}{3s+2}.
\end{equation}

This means that the analytic solution to Wiener index of generalized version $\mathcal{V}_{t,s}$ upon an arbitrary allowed seed is certainly obtained.

To make further progress, using Eqs.(\ref{eqa:MF-2-C-1}) and (\ref{eqa:MF-2-C-10}), the exact solution of mean first-passage time $\mathcal{A}_{\mathcal{V}_{t,s}}$ is shown in corollary 4.

\textbf{Corollary 4} The exact solution of mean first-passage time $\mathcal{A}_{\mathcal{V}_{t,s}}$ on generalized version $\mathcal{V}_{t,s}$ is

\begin{equation}\label{eqa:MF-2-C-11}
\mathcal{A}_{\mathcal{V}_{t,s}}=\frac{2\times3^{t}(s+1)^{t}}{n}\mathcal{W}+n(3^{t}-1)(s+1)^{t-2}+\frac{2[3^{t}(s+1)^{t}-1]}{(3s+2)(s+1)}.
\end{equation}

By far, we have accomplished the computation of mean first-passage time $\mathcal{A}_{\mathcal{V}_{t,s}}$ for random walks on generalized version $\mathcal{V}_{1,s}$.

As a case study, considering typical Vicsek fractal $\mathrm{V}_{t,s}$ with parameters $n=s+1$ and $\mathcal{W}=(s+1)^{2}$, we can reach a published result, i.e., the precise expression of quantity $\mathcal{A}_{\mathrm{V}_{t,s}}$.

\textbf{Corollary 5} The precise expression of quantity $\mathcal{A}_{\mathrm{V}_{t,s}}$ of typical Vicsek fractal $\mathrm{V}_{t,s}$ is expressed

\begin{equation}\label{eqa:MF-2-C-12}
\mathcal{A}_{\mathrm{V}_{t,s}}=2\times3^{t}(s+1)^{t+1}+(3^{t}-1)(s+1)^{t-1}+\frac{2[3^{t}(s+1)^{t}-1]}{(3s+2)(s+1)}.
\end{equation}
A more detailed calculation is deferred in Appendix.

\subsection{Scaling relation}
It is worth noticing that the closed form shown in Eq.(\ref{eqa:MF-2-C-12}) has been captured in \cite{Wu-2012} in terms of the self-similar property of Vicsek fractal $\mathrm{V}_{t,s}$. Yet, the generalized formula in Eq.(\ref{eqa:MF-2-C-11}) is not easily obtained using the methods introduced in \cite{Wu-2012} mainly because the seed used is not a star but an arbitrary tree $T$.

\begin{figure*}
\centering
\subfigure[]{
\begin{minipage}[t]{0.33\linewidth}
\centering
\includegraphics[width=6cm]{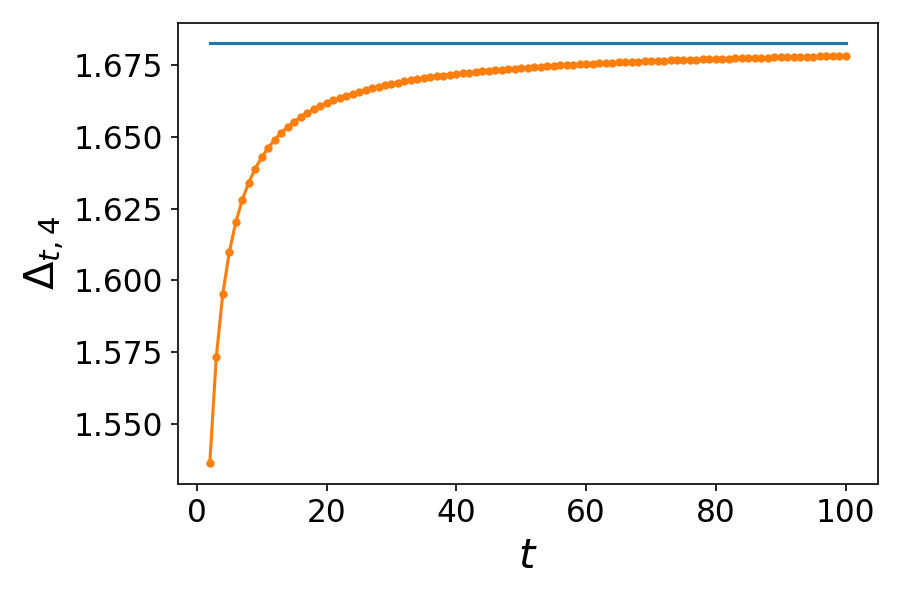}
%\caption{fig1}
\end{minipage}
}%
\subfigure[]{
\begin{minipage}[t]{0.33\linewidth}
\centering
\includegraphics[width=6cm]{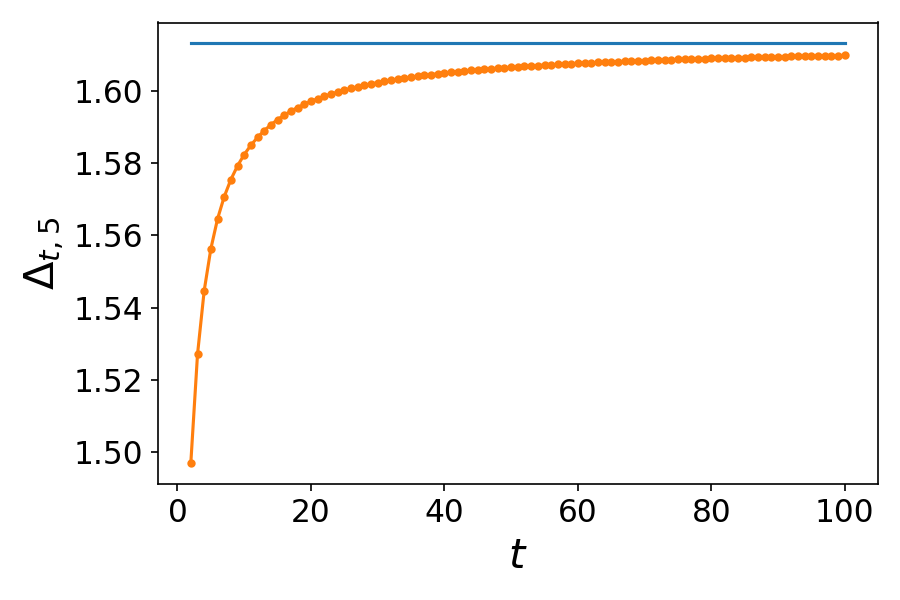}
%\caption{fig2}
\end{minipage}
}%
\subfigure[]{
\begin{minipage}[t]{0.33\linewidth}
\centering
\includegraphics[width=6cm]{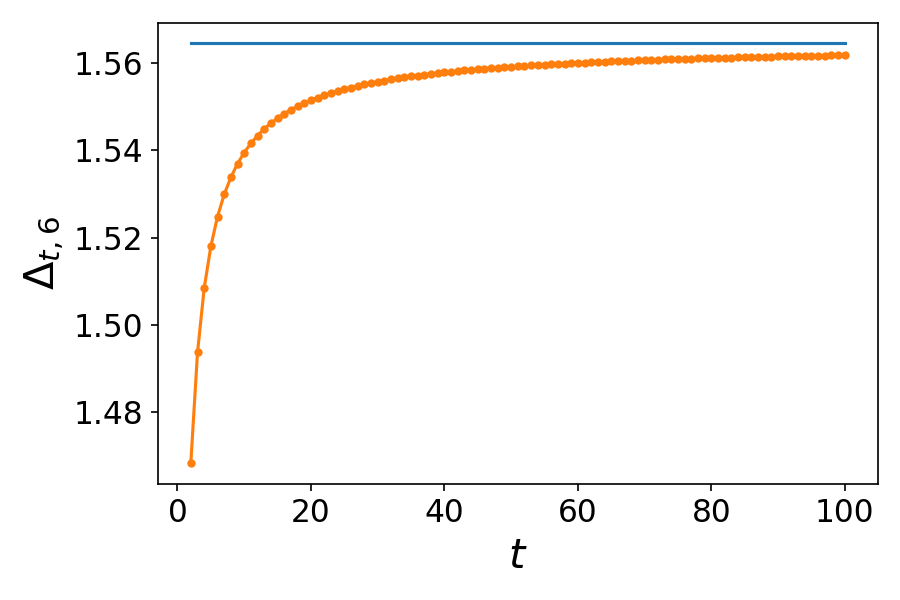}
%\caption{fig3}
\end{minipage}
}%
\\
\subfigure[]{
\begin{minipage}[t]{0.33\linewidth}
\centering
\includegraphics[width=6cm]{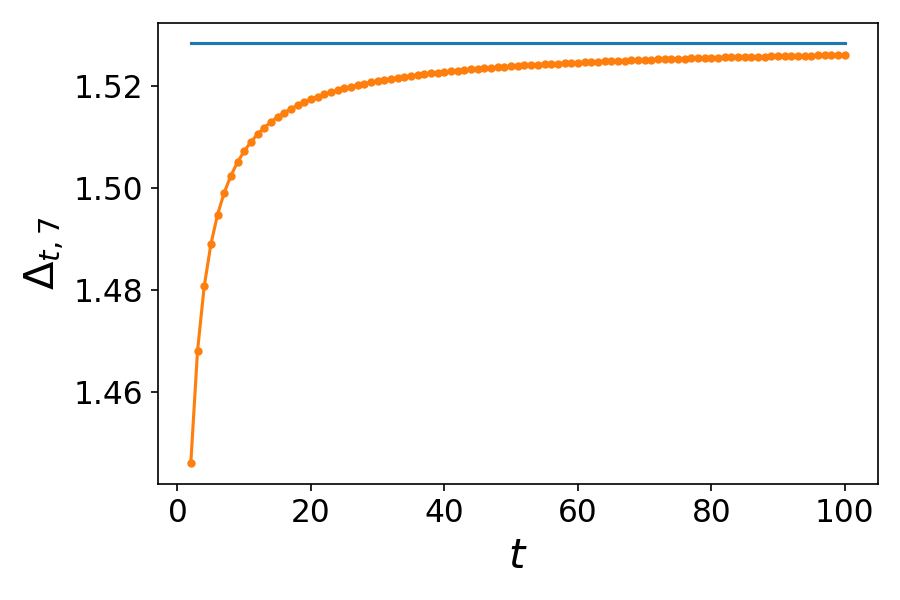}
%\caption{fig1}
\end{minipage}
}%
\subfigure[]{
\begin{minipage}[t]{0.33\linewidth}
\centering
\includegraphics[width=6cm]{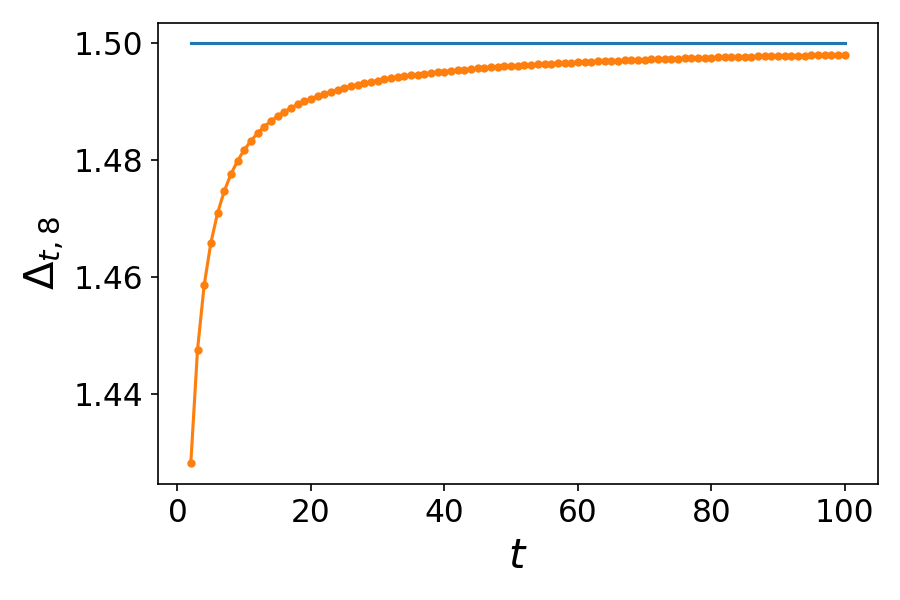}
%\caption{fig2}
\end{minipage}
}%
\subfigure[]{
\begin{minipage}[t]{0.33\linewidth}
\centering
\includegraphics[width=6cm]{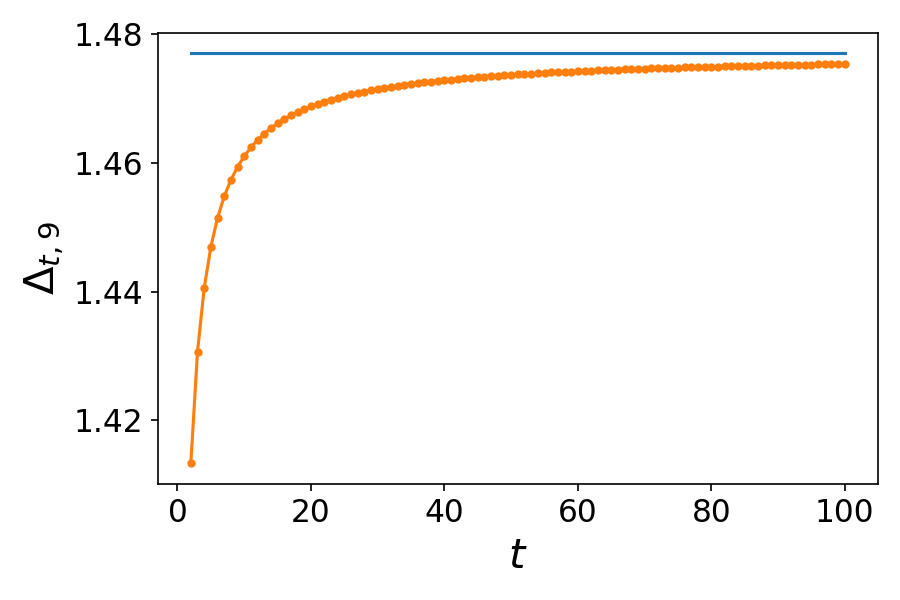}
%\caption{fig3}
\end{minipage}
}%
\caption{\label{fig:wide} (Color online)  The diagram of parameter $\Delta_{t,s}=\ln\mathcal{A}_{\mathcal{V}_{t,s}}/\ln|\mathcal{V}_{t,s}|$ where parameter $s$ is set to be $4,5,\cdots,9$, and we have made use of tree plotted in panel (a) of Fig.3 as a seed in which $n=5$ and $\mathcal{W}=18$ for purpose of illustration. The horizontal line indicates value for exponent $\lambda$ defined in Eq.(\ref{eqa:MF-2-C-14}). Obviously, the curve line representing value for parameter $\Delta_{t,s}$ gradually approaches the corresponding horizontal line in the large-$t$ limit when considering each value for parameter $s$, implying that experimental simulations are perfectly consistent with the theoretical analysis.}
\end{figure*}

In the limit of large graph size, the quantity $\mathcal{A}_{\mathcal{V}_{t,s}}$ in the preceding equation can be approximately equal to

\begin{equation}\label{eqa:MF-2-C-13}
\mathcal{A}_{\mathcal{V}_{t,s}}=O(3(s+1))^{t}.
\end{equation}

To make further process, substituting $\ln|\mathcal{V}_{t,s}|\sim \ln(s+1)\times t$ into Eq.(\ref{eqa:MF-2-C-13}) yields

\begin{equation}\label{eqa:MF-2-C-14}
\mathcal{A}_{\mathcal{V}_{t,s}}\sim|\mathcal{V}_{t,s}|^{\lambda}, \quad \lambda=\frac{\ln3}{\ln(s+1)}+1.
\end{equation}

Fig.4 illustrates the relation $\Delta_{t,s}=\frac{\ln\mathcal{A}_{\mathcal{V}_{t,s}}}{\ln|\mathcal{V}_{t,s}|}$ between quantities $\mathcal{A}_{\mathcal{V}_{t,s}}$ and $|\mathcal{V}_{t,s}|$ as a function of both the time step $t$ and parameter $s$, implying that the experimental simulations are in perfect agreement with the theoretical analysis.

In particular, when discussing the typical Vicsek fractal $\mathrm{V}_{t,s}$, the mean first-passage time $\mathcal{A}_{\mathrm{V}_{t,s}}$  can be asymptotically expressed in the following form

\begin{equation}\label{eqa:MF-2-C-15}
\mathcal{A}_{\mathrm{V}_{t,s}}\sim|\mathrm{V}_{t,s}|^{\lambda}, \quad \lambda=\frac{\ln3}{\ln(s+1)}+1.
\end{equation}

From both Eqs.(\ref{eqa:MF-2-C-14}) and (\ref{eqa:MF-2-C-15}), we can demonstrate that the choice of seed does not have a crucial effect on the exponent $\lambda$. Roughly speaking, the scaling parameter $\lambda$ is often closely related to the spectral dimension $d$ via formula $\lambda=2/d$. This suggests that the value for spectral dimension $d$ of generalized version $\mathcal{V}_{1,s}$ is calculated to write

\begin{equation}\label{eqa:MF-2-C-16}
d=\frac{2}{\lambda}=\frac{\ln(s+1)^{2}}{\ln3(s+1)}.
\end{equation}

Additionally, with another quality $d=2d_{f}/d_{w}$ where $d_{f}$ and $d_{w}$ represent the fractal dimension and random-walk dimension of generalized version $\mathcal{V}_{1,s}$, respectively, we may have

\begin{equation}\label{eqa:eqa:MF-2-C-17}
d_{w}=\frac{2d_{f}}{d}=1+\frac{\ln(s+1)}{\ln3}
\end{equation}
here the fractal dimension equal to $\ln(s+1)/\ln3$  has been used.

As noted in \cite{Gennes-1982}, the probability $P_{vv}(t)$ of returning back to vertex $v$, at long time $t$, follows the relation $P_{vv}(t)\sim t^{-d/2}$ when a walker starts out from a given vertex $v$ in graph $G(V,E)$ under consideration. So, in the large-$t$ limit, using parameter $d$ in Eq.(\ref{eqa:MF-2-C-16}) can allow us to better understand the distribution of returning probability. In addition, based on the same demonstration as in \cite{Gennes-1982}, we may find that a walker originating from a designed vertex will return back to the vertex almost surely over time when parameter $s$ is equivalent to $3$. This is because the corresponding spectral dimension $d$ in Eq.(\ref{eqa:MF-2-C-16}) is no more than $2$.

Last but not least, it is clear to understand that fractal usually exhibits self-similar property popularly observed in a lot of complex networks \cite{Song-2005,Zhang-2016}. As said above, the result in \cite{Wu-2012} is indeed established in terms of such a characteristic of the typical Vicsek fractal $\mathrm{V}_{t,s}$. As a result, some other exactly decidable fractals, such as T-graph, may also be studied in a similar manner proposed in this paper in order to estimate some interesting structural parameters, for example, Wiener index and mean first-passage time.

\section{Conclusion}

In conclusion, we revisit the typical Vicsek fractal $\mathrm{V}_{t,s}$, modeling a class of regular hyperbranched polymers, and study its corresponding generalized version $\mathcal{V}_{t,s}$ built upon an arbitrary allowed tree as a seed. Next, we consider random walks on model $\mathcal{V}_{t,s}$, and derive exact solution to mean first-passage time $\mathcal{A}_{\mathcal{V}_{t,s}}$. As opposed to previous works, we propose a mapping-based method to correctly evaluate quantity $\mathcal{A}_{\mathcal{V}_{t,s}}$ based on connection between Wiener index $\mathcal{W}$ and mean first-passage time $\mathcal{A}$ on tree. As a result, the pre-existing results related to Vicsek fractal $\mathrm{V}_{t,s}$ are covered by ours completely. More importantly, the method proposed is more convenient to carry out computation of quantity $\mathcal{A}_{\mathcal{V}_{t,s}}$ than many other methods including spectral techniques under a general situation. Note also that this method can be certainly employed to consider mean first-passage time on many other dendrimers, such as Cayley trees, for future work.

\section*{Acknowledgments}
We would like to thank Renbo Zhu for valuable assistance. The research was supported in part by the National Key Research and Development Plan under grant 2017YFB1200704 and the National Natural Science Foundation of China under grant No.61662066.

\section*{Appendix}

As a case study, we consider classic Vicsek fractal $\mathrm{V}_{t,s}$ and study its corresponding mean first-passage time $\mathcal{A}_{\mathrm{V}_{t,s}}$ using a similar manner used in \cite{Zhang-2010}.

Without loss of generality, let $\mathbf{L}_{t,s}$ denote by Laplacian matrix of Vicsek fractal $\mathrm{V}_{t,s}$ that is defined as

$$
l_{s,t}(ij)=\left\{\begin{aligned}&-1,\quad \text{vertex $i$ is connected to $j$ by an edge;}\\
&\quad k_{i},\quad \text{vertex $i$ is identical to $j$;}\\
&\quad0,\quad \;\:\text{otherwise.}
\end{aligned}
\right.
$$

As a result, its corresponding pseudoinverse, referred to as $\mathbf{L}_{t,s}^{\ddag}$, can be wrote in the following form

\begin{equation}\label{eqa:AP-0}
\mathbf{L}_{t,s}^{\ddag}=\left(\mathbf{L}_{t,s}-\frac{\mathbf{U}_{t,s}\mathbf{U}_{t,s}^{\top}}{|\mathrm{V}_{t,s}|}\right)^{-1}+\frac{\mathbf{U}_{t,s}\mathbf{U}_{t,s}^{\top}}{|\mathrm{V}_{t,s}|}.
\end{equation}
Here, symbol $\mathbf{U}_{t,s}$ is a $|\mathrm{V}_{t,s}|$-dimensional vector whose all entries are equivalent to $1$, i.e.,
$$\mathbf{U}_{t,s}=(\underbrace{1,1,...,1}_{|\mathrm{V}_{t,s}|})^{\top}.$$

Accordingly, the first-passage time $F_{u\rightarrow v}(s,t)$ for a couple of vertices $u$ and $v$ may be obtained in terms of entries $l^{\ddag}_{s,t}(ij)$ in matrix $\mathbf{L}_{t,s}^{\ddag}$, as below

\begin{equation}\label{eqa:AP-1}
F_{u\rightarrow v}(s,t)=\sum_{\alpha=1}^{|\mathrm{V}_{t,s}|}[l^{\ddag}_{s,t}(i\alpha)-l^{\ddag}_{s,t}(ij)+l^{\ddag}_{s,t}(j\alpha)-l^{\ddag}_{s,t}(jj)]l_{s,t}(\alpha\alpha)
\end{equation}
in which $l_{s,t}(\alpha\alpha)$ indicates the $\alpha$-th diagonal entry of Laplacian matrix $\mathbf{L}_{t,s}$.

Clearly, Vicsek fractal $\mathrm{V}_{t,s}$ is of tree structure and thus the exact solution to quantity $\mathcal{A}_{\mathrm{V}_{t,s}}$ can be expressed using all the non-zero eigenvalues $\varphi_{t,s;i}$ of its Laplacian matrix $\mathbf{L}_{t,s}$, as follows

\begin{equation}\label{eqa:AP-2}
\mathcal{A}_{\mathrm{V}_{t,s}}=\frac{\sum_{u,v;u\neq v}F_{u\rightarrow v}(s,t)}{|\mathrm{V}_{t,s}|(|\mathrm{V}_{t,s}|-1)}=2\sum_{i=2}^{|\mathrm{V}_{t,s}|}\frac{1}{\varphi_{t,s;i}}.
\end{equation}
This provides us with an approach to determining closed-form expression of $\mathcal{A}_{\mathrm{V}_{t,s}}$. Using some previous results based on real-space decimation method reported in \cite{Jayanthi-1993} and \cite{Jayanthi-1992}, we can find a recursive connection between two successive eigenvalues of an unique vertex $i$

\begin{equation}\label{eqa:AP-3}
\varphi_{t+1,s;i}(\varphi_{t+1,s;i}-3)(\varphi_{t+1,s;i}-s-1)=\varphi_{t,s;i}.
\end{equation}
Once the initial condition for eigenvalue $\varphi_{t,s;i}$ is not equal to zero, taking into account elementary rule in algebra, one can make use of necessary technique to solve for $\varphi_{t,s;i}$ from Eq.(\ref{eqa:AP-3}) and obtain three new and distinct values, $\varphi^{1}_{t,s;i}, \varphi^{2}_{t,s;i}$ and $\varphi^{3}_{t,s;i}$, at each time step $t$. More importantly, those newly generated eigenvalues keep the degeneracy of their ancestors. Along such a line, one might obtain all non-zero eigenvalues of its Laplacian matrix of Vicsek fractal $\mathrm{V}_{t,s}$. Nonetheless, it should be noted that it is not easy to derive the desired results about explicit formulas of eigenvalues in an iterative manner mentioned above.

In order to address this issue, Zhang \emph{et al} in \cite{Zhang-2010} divided all non-zero eigenvalues $\varphi_{t,s;i}$ into two disjoint sets, $A_{t,s}^{1}$ and $A_{t,s}^{2}$, with respect to the fact that a portion of these eigenvalues are  nondegenerate and the remainders are degenerate. The former $A_{t,s}^{1}$ consists of all nondegenerate eigenvalues. The other eigenvalues are collected in set $A_{t,s}^{2}$. Afterwards, Eq.(\ref{eqa:AP-2}) may be reorganized as

\begin{equation}\label{eqa:AP-4}
\mathcal{A}_{\mathrm{V}_{t,s}}=2\left(\sum_{i\in A_{t,s}^{1}}\frac{1}{\varphi_{t,s;i}}+\sum_{j\in A_{t,s}^{2}}\frac{1}{\varphi_{t,s;j}}\right).
\end{equation}

By using somewhat complicated methods (see \cite{Zhang-2010} for more details), one can obtained

\begin{equation}\label{eqa:AP-5}
\left\{\begin{aligned}&\sum_{i\in A_{t,s}^{1}}\frac{1}{\varphi_{t,s;i}}=\frac{(3s+3)^{t}-1}{(s+1)(3s+2)}\\
&\sum_{j\in A_{t,s}^{2}}\frac{1}{\varphi_{t,s;j}}=\frac{(s-2)(s+1)^{t-1}(3^{t}-1)}{2}+\frac{(3s+3)^{t}-1}{3s+2}
\end{aligned}
\right.
\end{equation}
Plugging results from Eq.(\ref{eqa:AP-5}) into Eq.(\ref{eqa:AP-4}) yields our desirable solution.

As previously, it is straightforward to see that the computation of quantity $\mathcal{A}_{\mathrm{V}_{t,s}}$ for typical Vicsek fractal upon eigenvalues is remarkably complicated than the mapping-based method proposed by us in this paper. More generally, this type of methods will be prohibitive when an arbitrary tree is selected as a seed for generating generalized Vicsek fractal $\mathcal{V}_{t,s}$. Fortunately, our methods are adequately exploited to address such issues in more general situation as shown above.

\end{document}